\documentstyle[12pt]{amsart}
\title[the Kaplansky conjecture]{Weak integral forms 
and the sixth Kaplansky conjecture}
\author{Dmitriy Rumynin}
\address{Mathematics Department, 
University of Warwick, Coventry, CV4 7AL, U.K.}
\email{D.Rumynin@@warwick.ac.uk}
\subjclass{Primary 16W30}
\date{May 21, 1998}
\thanks{The research was partially supported by NSF}
\newtheorem{theorem}{Theorem}

\newtheorem{lemma}[theorem]{Lemma}

\newtheorem{cor}[theorem]{Corollary}

\newcommand{\Q}{{\mathbb Q}}
\newcommand{\A}{{\mathbb A}}

\begin{document}

\begin{abstract}
  It is a short unpublished note from 1998. I make it public because Cuadra and Meir refer to it in their paper \cite{CM}.
  
We precisely state and prove a folklore result
that if a finite dimensional semisimple Hopf
algebra 
admits a weak integral form then
it is of Frobenius type.
We use an argument similar to that of Fossum \cite{fos},
which predates the Kaplansky conjectures.
\end{abstract}

\maketitle

Let $H$ be a  semisimple Hopf algebra of dimension $n$
over an algebraically closed field $K$ of zero characteristic.
We say that $H$ is of Frobenius type if the degree of
of any irreducible representation of $H$ divides $n$.
The sixth Kaplansky conjecture \cite{kap} is that no Hopf algebra
not of Frobenius type exists.

The Hopf algebra operations on $H$ are denoted $m$, $u$, $\Delta$,
$\varepsilon$ and $S$. 
Let us choose $\Lambda\in\int_H$ 
such that $\varepsilon (\Lambda)=n$.
We define
$\nu\in H\otimes_KH$ by 
$\nu=\sum_{(\Lambda)}\Lambda_1\otimes S(\Lambda_2)$. 
We fix a simple $H$-module $M$ of dimension $k$ from now on.
The representation is $\rho:H\rightarrow\mbox{End}_KM$. 
Its character is $\chi=\mbox{Tr}\circ\rho:H\rightarrow K$.
We denote the indecomposable central idempotent corresponding
to $M$ by $e$.

\begin{lemma}
The identity
$e=\frac{k}{n}\sum_{(\Lambda)}\chi(\Lambda_1)S(\Lambda_2)$
holds.
\end{lemma}

{\bf Proof.}
Let $z=\frac{k}{n}\sum_{(\Lambda)}\chi(\Lambda_1)S(\Lambda_2)$.
Using Larson's orthogonality relations \cite{lr1},
we see that
$$ 
\chi(z)=
\frac{k}{n}\sum_{(\Lambda)}\chi(\Lambda_1)\chi(S(\Lambda_2))=
k\langle \chi,\chi\rangle=k, 
$$
while for a different irreducible character $\eta$
$$
\eta(z)=
\frac{k}{n}\sum_{(\Lambda)}\chi(\Lambda_1)\eta(S(\Lambda_2))=
k\langle \chi,\eta\rangle=0.
$$
Thus, it suffices to show that $z$
belongs to the center of $H$. 
Pick $h\in H$. We recall that $S^2=\mbox{Id}_H$
by the Larson-Radford theorem \cite{lr2}. 
We also note that $\chi(ab)=\chi(ba)$
because of the similar trace property. 
We start with the equality
\begin{eqnarray*}
&\sum_{(h)}S(h_1)\Lambda h_2=\varepsilon (h)\Lambda;
& \mbox{ then apply } \Delta: \\
&\sum_{(h,\Lambda)}S(h_2)\Lambda_1 h_3\otimes S(h_1)\Lambda_2h_4
=\sum_{(\Lambda)}\varepsilon (h)\Lambda_1\otimes\Lambda_2;
& \mbox{ apply } \mbox{Id}\otimes S: \\
&\sum_{(h,\Lambda)}S(h_2)\Lambda_1h_3 \otimes S(h_4)S(\Lambda_2)h_1
=\sum_{(\Lambda)}\varepsilon (h)\Lambda_1\otimes S(\Lambda_2);
& \mbox{ apply } \chi\otimes\mbox{Id}: \\
&\sum_{(h,\Lambda)} \chi(\Lambda_1) \otimes S(h_2)S(\Lambda_2)h_1
=\sum_{(\Lambda)}\varepsilon (h)\chi(\Lambda_1) \otimes S(\Lambda_2);
& \mbox{ apply } m: \\
&\sum_{(h)} S(h_2)zh_1
= \varepsilon (h)z.
& 
\end{eqnarray*}
This implies that $z$ is central since
$$zh=\sum_{(h)}h_3S(h_2)zh_1= 
\sum_{(h)}h_2\varepsilon(h_1)z =hz. \ \ \Box$$

Let $R$ be a subring of $K$. 
By a weak $R$-form of $H$, we understand
an $R$-form of $(H,m,\nu)$. 
This is a free $R$-submodule $H_R\subseteq H$
with an $R$-basis $x_1, \,\ldots\,,x_n$
that is a $K$-basis of $H$ and 
the coefficients $m_{i,j}^k$ and $\nu^{i,j}$
belong to $R$. These coefficients appear if one writes down
$m$ and $\nu$ in terms of the basis,
i.e. $m(x_i\otimes x_j)=m_{i,j}^kx_k$ 
and $\nu=\nu^{i,j}x_i\otimes x_j$
keeping summation by the matching 
sub- and super-scripts in mind.
The proof of the following theorem is 
based on the idea of \cite{fos}.

\begin{theorem}
Assume that $R\subseteq K$ is a UFD (unique factorization
domain). If $H$ admits a weak $R$-form then $\frac{n}{k}\in R$.
\end{theorem}

{\bf Proof.}
Pick $x\in M$. Let $N=H_R\cdot x$. The $R$-module $N$ 
is torsion-free 
and, therefore, free since $R$ is a UFD. 
Let $e_1, \,\ldots\,,e_k$ be an $R$-basis of $N$.
Since $N$ is simple, it is also a $K$-basis of $M$.

We can write elements $\rho(t)$ with $t\in H$ as matrices
in this basis. If $t\in H_R$ then all coefficients of
$\rho(t)$ belong to $R$. 
In particular, $\chi(t)$ belongs to $R$.
%The integral $\Lambda$ belongs to $H_R$ since
%$\nu$ is defined over $R$. There exist $a_i,b_i\in H_R$
%such that $\sum_{(\Lambda)} \Lambda_1 \otimes S(\Lambda_2) =
%\sum_ia_i\otimes b_i$
%since $\Delta$ and $S$ are defined over $R$. 
By Lemma 1,   
$$\frac{n}{k}\mbox{Id}_V= \frac{n}{k}\rho(e) = 
\chi\otimes\rho(\nu)=\nu^{i,j}\chi(x_i)\rho(x_j)$$
The right part is apriori a matrix with coefficients
in $R$. Thus, $\frac{n}{k}\in R$ 
as a coefficient of the left part.
$\ \ \Box$

Let $\A\subseteq \overline{\Q}\subseteq K$ be the ring 
of integer algebraic
numbers. The following corollary is straightforward. 

\begin{cor}
If there exists a collection
$R_i$ of UFDs such that $H$ admits a weak $R_i$-form for each $i$
and $\cap_iR_i\subseteq \A$ then
$H$ is of Frobenius type.
\end{cor}

\begin{cor} If $H$ admits a weak $\A\cap L$-form
for an algebraic number field $L$ then $H$
is of Frobenius type.
\end{cor}

{\bf Proof.} It suffices to consider the collection
of discrete valuation rings of $L$ for all $p$-adic
valuations and to apply Corollary 3. $\ \Box$

\end{document}